\documentclass[12pt]{article}
\usepackage{a4}
\usepackage{amssymb}
\usepackage{amsmath}
\usepackage{amsthm}
\usepackage{amstext}
\usepackage{amscd}
\usepackage{latexsym}
\usepackage{graphics}
\theoremstyle{plain}
\newtheorem{thm}{Theorem}[section]
\newtheorem{lemma}[thm]{Lemma}

\newtheorem{cor}[thm]{Corollary}
\newtheorem{prop}[thm]{Proposition}
\theoremstyle{definition}
\newtheorem{rmk}[thm]{Remark}

\newtheorem{notation}[thm]{Notation}
\newtheorem{ex}[thm]{Example}

\newcommand{\tensor}{\otimes}

\newcommand{\Union}{\bigcup}

\newcommand{\disc}{\partial}

\renewcommand{\tilde}{\widetilde}
\newcommand{\p}{{\mathfrak p}}
\newcommand{\q}{{\mathfrak q}}

\newcommand{\sO}{{\mathcal O}}

\newcommand{\F}{{\mathbb F}}

\newcommand{\N}{{\mathbb N}}

\newcommand{\Q}{{\mathbb Q}}

\newcommand{\Z}{{\mathbb Z}}
\newcommand{\Gal}{{\rm Gal}}
\input{xy}
\xyoption{all}
\begin{document}

\title{Absolute norms of $p$-primary units}
\author{Supriya Pisolkar}
\date{}
\maketitle

\section{Introduction}
 Let $\rm{K}$ be a $p$-adic field containing a primitive $p$-th root of unity $\zeta_p$. A unit $\alpha \in \rm{U_K}$ is called {\bf $p$-primary} if the extension ${\rm K}(\alpha^{1/p})$ is an unramified extension of ${\rm K}$. For example, the discriminant of an integral basis of an unramified extension $\rm{L/K}$ of a $2$-adic field is a $2$-primary unit since it has a square root in $\rm{L}$ (see proof of Cor.\ref{p=2}). It is interesting to observe that for an extension $\rm{ L/K}$ the norm of a $p$-primary unit in $\rm L$ is again a $p$-primary unit in $\rm K$ (see Lemma \ref{primary-norm}). In this paper we prove a result about the absolute norm of a $p$-primary unit which is motivated by the following theorem of J. Martinet. 
\begin{thm} \label{martinettheorem} {\rm(\cite{martinet}; 1.4)} Let $\rm{L/K}$ be a quadratic extension of number fields such that the absolute norm of the relative discriminant ideal $\disc_{\rm L/K}$ is odd. If $\rm K$ contains the $2^m$-th roots of unity for some $m\geq 2$, then  ${\rm N}_{\rm K/\Q}(\disc_{{\rm L/K}}) \equiv 1({\rm mod} \ 2^{m+1})$.
\end{thm}

\noindent Although the above result has been stated only for a quadratic extension, arguments in \cite{martinet} are sufficient to obtain the same congruence for any finite extension of number fields satisfying the above hypothesis, by reducing to the quadratic case. The above result is about the behaviour of the absolute norm of the discriminant ideal at the prime $2$. Thus one may ask if there is a local analogue of Martinet's theorem.   \\  

\noindent The main contribution of this paper is to prove a local analogue (see corollary \ref{p=2}) of Theorem \ref{martinettheorem} which turns out to be a statement about $2$-primary units. In fact, we prove the following theorem about $p$-primary units for all primes $p$.

\begin{thm}\label{main}
Let $\rm{K}$ be a finite extension of $\Q_p$ containing a primitive $p^m$-th root of unity where $m\geq 1$. Let $\alpha$ be a $p$-primary unit of $\rm {U_K}$. Then $${\rm N}_{{\rm K}/\Q_p}(\alpha) \equiv 1({\rm mod} \ {\it p}^{m+1})$$
\end{thm}
\begin{rmk}
The absolute norm of a $p$-primary unit may not satisfy a better congruence. Moreover, this congruence may not hold for arbitrary units. See example \ref{examples}.
\end{rmk}

\begin{cor}\label{p=2} Let $\rm{K}$ be a finite extension of $\Q_2$ containing a primitive $2^m$-th root of unity for some $m \geq 1$. Let ${\rm L}$ be an unramified extension of $\rm{K}$ and let $\{\alpha_1, \alpha_2, \ldots, \alpha_n\}$ be an integral basis of ${\rm L/K}$.Then, 
$${\rm N}_{{\rm K}/\Q_2}(d_{\rm L/K}(\alpha_1,\alpha_2, \ldots ,\alpha_n)) \equiv 1({\rm mod} \ 2^{m+1})$$  
\end{cor}

\begin{proof} $\ {\rm K}(\sqrt{d_{\rm L/K}(\alpha_1,\ldots,\alpha_n)}\ ) \ $ is a subextension of ${\rm L/K}$ and thus an unramified extension of ${\rm K}$. Therefore $d_{\rm L/K}(\alpha_1,\ldots,\alpha_n)$ is a $2$-primary unit of ${\rm K}$ and the result follows from the above theorem. 
\end{proof}

\noindent Using the above local statement, one would now like to recover Theorem \ref{martinettheorem}. We are able to do this in the case when ${\rm L/K}$ has an integral basis (see corollary \ref{deduce}). We first prove the following. 
\begin{cor}\label{global}
Let ${\rm K}$ be a number field containing the $2^m$-th roots of unity where $m\geq 1$. Let ${\rm L/K}$ be a finite extension, let $\disc_{\rm L/K}$ be the relative discriminant of $\rm L/K$, and suppose that the absolute norm ${\rm N_{K/\Q}}(\disc_{\rm L/K})$ is odd. Assume that $\underline{\alpha}= \{\alpha_1,...,\alpha_n\}$ is a basis of $\rm L/K$ such that
\begin{enumerate}
\item[(i)] $\underline{\alpha}$ is contained in $\sO_{\rm L}$. 
\item[(ii)] $\underline{\alpha}$ forms a $S^{-1}\sO_{\rm K}$-basis of $ S^{-1}{\sO}_{\rm L}$ where $S = \Z \backslash(2)$.
\end{enumerate}
\noindent Let $d_{\rm L/K}(\underline{\alpha})$ denote the discriminant of the basis $\underline{\alpha}$. Then 
$${\rm N_{K/\Q}}(d_{\rm L/K}(\underline{\alpha})) \equiv 1 ({\rm mod} \ 2^{m+1})$$
\end{cor}

\begin{proof}[Proof of Corollary \ref{global}] 
The hypothesis that the norm of the discriminant ideal is odd is equivalent to saying that all primes of $\sO_{\rm K}$ lying above the prime ideal $(2)$ are unramified in ${\rm L}$. Fix a prime $\p$ of $\sO_{\rm K}$ lying above $2$. Let $\{\q_j\}_{j=1}^r$ be the prime ideals of $\sO_{\rm L}$ lying above $\p$. Choose an integral basis $\underline{\beta_j}= \{ \beta_{j1},\beta_{j2},...\}$ of ${\rm L}_{j}/{\rm K}_{\p}$. Here ${\rm L}_j$ and ${\rm K}_{\p}$ denote the completion of ${\rm L}$ and ${\rm K}$ at $\q_j$ and $\p$ respectively. Let $A={\rm L} \tensor_{\rm K}{\rm K}_{\p}$. Since $A=\oplus{\rm L}_j$, $\underline{\beta}= \displaystyle\Union_{j=1}^{r} \underline{\beta_j}$ forms an integral basis for the ${\rm K}_{\p}$ algebra $A$. Here by integral basis of $A$ we mean an $\sO_{\rm K_{\p}}$-basis of $\oplus_{j} \sO_{{\rm L}_j}$. Let $d_{A/{\rm K}_{\p}}(\underline{\beta})$ denote the discriminant of $A$ of $\underline{\beta}$. Note that $\underline{\alpha}$ is also an integral basis of $A$ and the $\p$-adic image $(d_{\rm L/K}(\underline{\alpha}))_{\p}$ is equal to $d_{A/{\rm K}_p}(\underline \alpha)$. Therefore 
$$(d_{\rm L/K}(\underline{\alpha}))_{\p}= d_{A/{\rm K}_{\p}}(\underline{\beta})\cdot u^2$$ 
for some unit $u \in {\rm U}_{\rm K_{\p}}$. Since $d_{A/{\rm K}_{\p}}(\underline{\beta})= \displaystyle\prod_{j=1}^r d_{{\rm L}_j/{\rm K}_{\p}}(\underline{\beta_j})$, by Corollary \ref{p=2} we have, $${\rm N_{K_{\p}/\Q_2}}d_{A/{\rm K}_{\p}}(\underline{\beta}) =\displaystyle\prod_{j=1}^r {\rm N}_{ {\rm K}_{\p}/\Q_2}(d_{{\rm L}_j/{\rm K}_{\p}}(\underline{\beta_j})) \equiv 1 ({\rm mod}\ 2^{m+1})$$ 
As $\Q_2(\zeta_{2^m})\subset {\rm K}_{\p}$, by theorem \ref{zetam} we have, ${\rm N}_{{\rm K}_{\p}/\Q_2}(u) \in {\rm U}_{m,\Q_2}$. This implies that $ {\rm N}_{{\rm K}_{\p}/\Q_2}(u)^2 \in {\rm U}_{m+1,\Q_2}$. Thus $$ {\rm N_{K_{\p}/\Q_2}}(d_{\rm L/K}(\underline{\alpha}))_{\p} \equiv 1 ({\rm mod}\ 2^{m+1})$$
 By using (\cite{serre},II,3.2) we get $$ {\rm N_{K/\Q}}(d_{\rm L/K}(\underline{\alpha})) = \displaystyle\prod_{\p |2 }{\rm N_{K_{\p}/\Q_2}}(d_{\rm L/K}(\underline{\alpha}))_{\p} \equiv 1 ({\rm mod}\ 2^{m+1})$$
\end{proof}

\begin{cor} \label{deduce} Let ${\rm K}$ be as above with $m \geq 2$. Let $\rm L/K$ be a finite extension having an integral basis $\underline{\alpha}$. Assume that the absolute norm ${\rm N_{K/\Q}}$ of the discriminant ideal $\disc_{\rm L/K}$ is odd. Then ${\rm N_{K/\Q}}(d_{\rm L/K}(\underline{\alpha}))$ is the positive generator of ${\rm N_{K/\Q}}(\disc_{\rm L/K})$. In particular  $${\rm N_{K/\Q}}(\disc_{\rm L/K}) \equiv 1 ({\rm mod} \ 2^{m+1})$$
\end{cor}
 
\begin{proof} Since $\rm N_{{\rm K}/\Q}(d_{\rm L/K}(\underline \alpha))$ is a generator of ${\rm N_{K/\Q}}(\disc_{\rm L/K})$, it suffices to show that it is positive.  As $\rm K$ is purely imaginary, the norm $\rm {N}_{{\rm K}/\Q}$ of any nonzero element is positive.  
\end{proof}

\noindent Theorem \ref{main} will be proved in section \ref{sec-main}. In section 2, we recall some classical results on the behaviour of the norm map of a totally ramified cyclic extension of local fields. See for example \cite{fesenko}. \\

\noindent {\bf Acknowledgement}\  : \ I am thankful to Prof. C. S. Dalawat for his help. This problem was suggested to me by him. I am indebted to Jo\"{e}l Riou for painstakingly going through the first draft of this paper and for important suggestions. I am grateful to Amit Hogadi for his interest and stimulating discussions.  I thank Prof. Lo\"{i}c Merel for his useful comments.


\section{Cyclic ramified extensions and the norm map.} \label{results} 
In this section we state some results about the norm map of a totally ramified cyclic extension of local fields. These results will be crucially used in the proof of the theorem \ref{main}. 

\begin{notation} Let $\rm{K}$ be a finite extension of $\Q_p$. Let $k$ be the residue field of $\rm{K}$. We denote by ${\rm U}_{i,{\rm K}}$ the subgroup of ${\rm U}_{\rm K}$ given by 
$$ {\rm U}_{i,{\rm K}}= \{ x \in {\rm U_K} | v_{\rm K}(1-x) \geq i \}$$ 
We denote by $e_{\rm K}$ the absolute ramification index of ${\rm K}$ and ${\overline e}_{\rm K} = e_{\rm K}/{p-1}$.\\

\begin{thm}{\rm (\cite{serre}, p.~212)}\label{p-power} Let $\rm K$ be a finite extension of $\Q_p$. For all $n > \overline{e}_{\rm K}$,
the map $( \ )^p: {\rm U}_n \rightarrow  {\rm U}_{n+e_{\rm K}}$ is a bijection.
\end{thm}

\end{notation}

\begin{thm}{\rm(\cite{fesenko}, III.1.4) } \label{fesenko1}
 Let $\rm{L/K}$ be a totally ramified Galois extension of degree $p$. Let $\pi_{\rm L}$ be a uniformiser of ${\rm L}$. Let $\sigma$ be a generator of $\Gal(\rm{L/K})$. 
Then $\sigma(\pi_{\rm L})/{\pi_{\rm L}} \in  {\rm U_{1,L}}$. 
Further, if $s$ is the largest integer such that 
$\sigma(\pi_{\rm L})/{\pi_{\rm L}} \in {\rm U}_{s,\rm{L}}$,  
then $s$ is independent of the uniformiser $\pi_{\rm L}$. 
\end{thm}

\noindent Thus the integer $s$, defined above depends only on the extension $\rm{L/K}$. Therefore we will denote it by $s({\rm L/K})$. Note that $s({\rm L/K})$ is the unique ramification break of ${\rm Gal}({\rm L/K})$.

\begin{thm}{ \rm(\cite{fesenko},III.2.3) }\label{fesenko2}
 $s({\rm L/K}) \leq p{\overline e}_{\rm L}$.
\end{thm}

\begin{ex}
  Let ${\rm K}=\Q_2$. Then $\rm K$ has six ramified quadratic extensions namely, ${ \rm K(\sqrt{-1}), K(\sqrt{-5}), K(\sqrt{\pm 2}), K(\sqrt{\pm 10})} $. For the first two extensions $s=1$ and for the remaining extensions $s=2$.  
\end{ex}
\noindent We will now state some results about the Hasse-Herbrand function, which is an important tool in understanding the behaviour of the norm map in wildly ramified extension. 

\begin{prop}{\rm (\cite{fesenko}, III, prop.~3.1)} Let ${\rm L/F}$ be a finite Galois extension of local fields and let the residue field ${\rm k_F}$ to be infinite , ${\rm N} = {\rm N_{L/F}}$. The there exists a unique function  
$$ h = h_{\rm L/F}: \N \to \N$$  
such that $h(0)=0$ and 
$${\rm NU}_{h(i),\rm L} \subset {\rm U}_{i,{\rm F}},\ {\rm NU}_{h(i),\rm L} \not\subset {\rm U}_{i+1,\rm F}, \ {\rm NU}_{h(i)+1, \rm L} \subset {\rm U}_{i+1, \rm F}.$$
\end{prop}

\begin{enumerate}
\item[(1)] For ${\rm L/F}$ a totally tamely ramified extension, 
$$h(i) = [{\rm L}:{\rm F}]i,$$

\item[(2)] If ${\rm L/F}$ is totally ramified extension of degree $p = \text{char}({\rm k_F})$ then,
 
 \[ h(i) = \left\{ \begin{array}{ll} i, & \mbox{if $i < s$ }\\
                                s(1-p)+pi & \mbox{if $i\geq s$}, \end{array} \right. \]  

\end{enumerate}
\noindent If we have the tower of field extension ${\rm L} \subset {\rm M} \subset {\rm F}$ then, $$ h_{\rm L/F} = h_{\rm L/M} \circ h_{\rm M/F}.$$

\vspace{.1cm}
\noindent To treat the case of local fields with the finite residue fields we have the following.
\begin{lemma}{\rm (\cite{fesenko}, III,~3.2)} Let ${\rm L/F}$ be a finite separable totally ramified extension of local fields. Then for an element $\alpha \in {\rm L}$ we get 
$$ {\rm N_{L/F}}(\alpha) = {\rm N}_{\widehat{{\rm L}^{ur}}/\widehat{{\rm F}^{ur}}}.$$
where $\widehat{{\rm F}^{ur}}$ is the completion of ${\rm F}^{ur}$, $\widehat{{\rm L}^{ur}} = {\rm L}\widehat{{\rm F}^{ur}}$.
\end{lemma}

\noindent This lemma shows that, if ${\rm k_F}$ is finite then, for a finite Galois extension ${\rm L/F}$, $h_{\rm L/F} = h_{\widehat{{\rm L}^{ur}}/\widehat{{\rm L}^{ur}}}$. 

\noindent We will be using the following theorem from class field theory, 
\begin{thm}{\rm (\cite{neukirch2}, p.~45)}\label{zetam} The norm map carries units of $\Q_p(\zeta_{p^m})$ into ${\rm U}_{m,\Q_p}$; i.e, ${\rm N}_{\Q_p(\zeta_{p^m})\mid\Q_p}\left({\rm U}_{\Q_p(\zeta_{p^m})}\right) = {\rm U}_{m,\Q_p}$. 
\end{thm}

\noindent Next we recall some results discussed in \cite{dalawat}. Let $\rm K$ be a finite extension of $\Q_p$. We know that $\rm K^*$ has a filtration $({\rm U}_n)_{n \in \N}$. Put ${\rm U}_0 = {\rm U}_{\rm K}$ and $\overline{{\rm U}}_n = {\rm U}_n/{\rm U}_n \cap {\rm U}_0^p$. We thus a get a filtration on ${\rm K}^*/{\rm K^*}^p$ by $\F_p$ subspaces,
$$\cdots \subset {\overline{\rm U}}_n \subset \cdots \subset {\overline{\rm U}}_1 \subset {\overline{\rm U}}_0 \subset {\rm K^*}/{\rm K^*}^p$$ 

\noindent From \rm(\cite{dalawat}, Prop. 33), if $\rm K^*$ contains an element of order $p$ then $\overline{\rm U}_{p{\overline e}_{\rm k}}$ is an $\F_p$-line in ${\rm K}^*/{\rm K^*}^p$. Further, $\F_p$-lines in ${\rm K}^*/{\rm K^*}^p$ are in bijection with the cyclic degree $p$ extensions of $\rm K$. The line which corresponds to the unramified extension is given by following proposition.

\begin{prop}{\rm (\cite{dalawat}, Prop. 16)} \label{fp-line} The $\F_p$-line in ${\rm K^*/{K^*}}^p$ which gives the unramified $(\Z/p\Z)$-extension of ${\rm K}$ upon adjoining p-th roots is ${\overline {\rm U}}_{p{\overline e_{\rm K}}}$.
\end{prop}

\noindent We thus get the following important corollary.
\begin{cor}\label{cor:fpline} A unit $\alpha$ in ${\rm U_{K}}$ is $p$-primary if and only if the image of $\alpha$ in ${\rm K}^*/{\rm K^*}^p$  belongs to $\overline{\rm U}_{p\overline{e}_{\rm K}}$. 
\end{cor}

\noindent  In the earlier version of this paper, the statement of the following lemma was hidden in the proof of Theorem \ref{main}. I thank Jo\"el Riou for his suggestion of stating it as a separate lemma and also providing an elegant proof.
\begin{lemma}\label{primary-norm} Let $\rm{E/F}$ be an extension of $p$-adic fields containing a primitive $p$-th root of unity. Then the norm map ${\rm N_{E/F}}:{\rm E}^*\to {\rm F}^*$ takes a $p$-primary unit to a $p$-primary unit.
\end{lemma}   
\begin{proof} It suffices to prove this result in two cases: $\rm E/F$ is a totally ramified extension and $\rm E/F$ is an unramified extension. Let $\alpha \in \rm{E}$ be a $p$-primary unit. By definition of a $p$-primary unit, there exists an unramified extension $\rm{E'/E}$ and $\beta \in \rm E'$ such that $\alpha = \beta^p$.\\
\noindent Suppose that ${\rm E/F}$ is a totally ramified extension. Then the extension ${\rm E'/E}$ corresponds to the extension of residue field of $\rm{E}$. The residue field of $\rm E$ is the same as that of $\rm F$. Thus there exists an unramified extension $\rm{F'/F}$ such that ${\rm E'}={\rm F'E}$. It is easy to see that ${\rm N_{E/F}}(\alpha) = {\rm N_{E'/F'}}(\alpha)$ and thus ${\rm N_{E/F}}(\alpha) = {\rm N_{E'/F'}}(\beta)^p$.\\
\noindent Suppose now ${\rm E/F}$ is an unramified extension. ${\rm E'/E}$ is also an unramified extension. Since these extensions are Galois, ${\rm N_{E/F}}(\alpha)$ is a product of all $\sigma(\alpha)$ where $\sigma \in {\rm Gal}({\rm E/F})$. For each $\sigma \in {\rm Gal}(\rm E/F)$, we choose $\tilde{\sigma} \in {\rm Gal}({\rm E'/F})$ which extends $\sigma$. Then ${\rm N_{E/F}}(\alpha)$ is a $p$-th power of a product $\tilde{\sigma}(\beta)$. This product lies in an unramified extension of ${\rm F}$ namely ${\rm E'}$ and thus ${\rm N_{E/F}}(\alpha)$ is a $p$-primary unit. 
\end{proof}

\section{Proof of Theorem \ref{main}}\label{sec-main}
\begin{proof}[Proof of Theorem \ref{main} ] By hypothesis, $\rm K$ contains the $p^m$-th roots of unity where $m \geq 1$. By above Lemma \ref{primary-norm}, the norm ${\rm N}_{{\rm K}/\Q_p(\zeta_{p^m})}$ of a $p$-primary unit in ${\rm K}$ is a $p$-primary unit in $\Q_p(\zeta_{p^m})$. Thus it is enough to prove the result in the special case ${\rm K} = \Q_p(\zeta_{p^m})$.  \\

\noindent Suppose that ${\rm K} = \Q_p(\zeta_{p^m})$. By corollary \ref{cor:fpline}, a unit $\alpha \in {\rm K}$ is $p$-primary if and only if $$\ \alpha=u\cdot w^p$$ where $u  \in {\rm U}_{p \overline{e}_{\rm K},{\rm K}}$ and $w \in {\rm U}_{\rm K}$. By theorem \ref{zetam}, ${\rm N}_{{\rm K}/\Q_p}(u)$ and ${\rm N}_{{\rm K}/\Q_p}(w)$ belong to ${\rm U}_{m,\Q_P}$.   
Thus by theorem \ref{p-power}, ${\rm N}_{{\rm K}/\Q_p}(w)^p \in {\rm U}_{m+1,\Q_p}$. To prove the theorem it now remains to show that  ${\rm N}_{{\rm K}/\Q_p}(u) \in {\rm U}_{m+1,\Q_p}$. We are going to show this by using the Hasse-Herbrand function. In fact, we will show that $h_{{\rm K}/\Q_p}(m) = p^m-1$. Then, by using the property of Hasse-Herbrand function we will get that, ${\rm N}_{{\rm K}/\Q_p}({\rm U}_{h(m)+1, {\rm K}}) = {\rm N}_{{\rm K}/\Q_p}({\rm U}_{p^m,{\rm K}})  \subset {\rm U}_{m+1, \Q_p}$.\\

\noindent Consider the tower of field extensions ${\rm K}={\rm K}_m \supset  {\rm K}_{m-1} \supset \cdots \supset {\rm K}_1$ where $ {\rm K}_i = \Q_p(\zeta_{p^i})$. Note that for each $2 \leq i \leq m$, ${\rm K}_{i}/{\rm K}_{i-1}$ is a wildly ramified cyclic extension of degree $p$, and ${\rm K}_1/\Q_p$ is tamely ramified cyclic extension of degree $p-1$. Let $v_{{\rm K}_i}$ be the surjective valuation of ${\rm K}_i$.  By {\rm(\cite{serre}, IV, Lemma~1(c))}, $s_{{\rm K}_i/{\rm K}_{i-1}} = p^{i-1}-1$. Indeed,  $v_{{\rm K}_i}(\sigma(\zeta_{P^i})-\zeta_{p^i}) = v_{{\rm K}_i}(\zeta_p\zeta_{p^i}-\zeta_{p^i}) = v_{{\rm K}_i}(\zeta_p-1)= p^{i-1}$.\\

\noindent Henceforth for simplicity of notation we write $s_i$ for $s({\rm K}_i/{\rm K}_{i-1})$.

\noindent {\bf Step~1:} For $m=1$, we want to prove that $h_{\rm K_1/K}(1) = p-1$. Since ${\rm K_1/K}$ is a tamely ramified extension of the degree $p-1$, this follows by the formula of Hasse-Herbrand function for tamely ramified extensions. See (1).

\noindent {\bf Step~2:} For $m \geq 2$, we want to show that $h_{{\rm K}_m/\Q_p}(m) = {p^m}-1$. We will use the transitivity of the Hasse-Harbrand function thorough the tower of field extensions. For simplicity of notation, we will write
 $h_{{\rm K}_1/\Q_p } = h_1$, $h_{{\rm K}_i/{\rm K}_{i-1}} = h_i$, and $h = h_m \circ \cdots \circ h_1$.
\begin{enumerate}
\item[(i)] We know from (1) that $h_1(m) = m(p-1)$. To compute $h_2(m(p-1))$,  first observe that $s_2 = p-1$ and $m(p-1) > s_2$. Now applying formula (2), 
$$h_2(m(p-1)) = (p-1)(1-p)+p(m(p-1)) = (m-1)p^2-(m-2)p-1.$$
\item[(ii)] For $1\leq n < m$, let us assume that $$ (*) \ \ \  h_n \circ h_{n-1} \circ \cdots \circ h_1(m) = (m-(n-1))p^n-(m-n)p^{n-1}-1$$ 
\noindent It is easy to check that $$(m-(n-1))p^n-(m-n)p^{n-1}-1 \geq  s_{n+1} = p^n-1.$$
\noindent Now we can apply the formula (2).
\begin{align*}
h_{n+1}  \circ (*)  & \ \  = \ \ h_{n+1}((m-(n-1))p^n-(m-n)p^{n-1}-1) \\
                     & \ \  = \ \ (p^n-1)(1-p)+p[(m-(n-1))p^n-(m-n)p^{n-1}-1]\\                     & \ \  = \ \ (m-n)p^{n+1}-(m-(n+1))p^n-1   
\end{align*}

\noindent Thus, for $n+1=m$, we get $h_m \circ \cdots \circ h_1(m) = h(m) = p^m-1$.  
\end{enumerate}

\noindent This proves that ${\rm N}_{{\rm K}/\Q_p}({\rm U}_{h(m)+1, \rm K}) = {\rm N}_{{\rm K}/\Q_p}({\rm U}_{p^m,\rm K})  \subset {\rm U}_{m+1, \Q_p}$ and thus $${\rm N}_{{\rm K}/\Q_p}(u) \in {\rm U}_{m+1,Q_p}.$$ This completes the proof of the Theorem \ref{main}.

\end{proof}

\begin{rmk} The result ${\rm N_{K/\Q_p}}({\rm U}_{p^m,\rm K})  \subset {\rm U}_{m+1, \Q_p}$ can also be derived without the explicit use of the function $h_{{\rm K}/\Q_p}$. This can be achieved by the repeated application of the fact {\rm (\cite{fesenko}, III,~1.5)} that, for a totally ramified cyclic degree-$p$ extension ${\rm L/K}$, we have 

\begin{enumerate}
\item ${\rm N_{L/K}}({\rm U}_{s+pi,{\rm L}}) \subseteq {\rm U}_{s+i,{\rm K}} \ \ \forall \ i>0$
\item ${\rm N_{L/K}}({\rm U}_{s+i,{\rm L}})={\rm N_{L/K}}({\rm U}_{s+i+1,{\rm L}})$ for $i>0$, $p \nmid i$.
\end{enumerate}
\noindent and for totally tamely ramified Galois extension of degree-$n$ {\rm (\cite{fesenko}, III,~1.3)},  
\begin{enumerate}
\item ${\rm N_{L/K}}({\rm U}_{ni,{\rm L}}) \subseteq {\rm U}_{i,{\rm K}}$
\item ${\rm N_{L/K}}({\rm U}_{i,{\rm L}})= {\rm N_{L/K}}({\rm U}_{i+1,{\rm L}})$ if $\ n\nmid i$.
\end{enumerate}
\noindent In fact, in the notation ${\rm K}_m = \Q_p(\zeta_{p^m})$, these formulae imply 
$${\rm N}_{{\rm K}_i/{\rm K}_{i-1}} ({\rm U}_{r_i,{\rm  K}_i}) \subset  {\rm U}_{r_{i-1},{\rm  K}_{i-1}}, \ \text{for} \  2 \leq i \leq m. $$
\noindent where for $1 \leq i \leq m$, $r_i = (m-i+1)p^i -(m-i)p^{i-1}$. It is easy to see that $${\rm N}_{{\rm K}_1/\Q_p}({\rm U}_{r_1,{\rm K}_1}) \subset {\rm U}_{m+1,\Q_p}.$$ 

\end{rmk}


\begin{ex}\label{examples}
Let us show that a $p$-primary unit may not satisfy a congruence better than that of Theorem\ref{main}. Suppose that ${\rm K} = \Q_3(\zeta_{3^4})$. Let $\alpha \in {\rm U}_{2,\Q_3} \backslash {\rm U}_{3,\Q_3}$. Then $\alpha \in {\rm U}_{3^4,{\rm K}}$ i.e $\alpha \in {\rm U}_{3\overline{e}_{\rm K}, {\rm K}}$. Therefore $\alpha$ is a 3-primary unit by using \ref{fp-line}. We claim that  
${\rm N}_{{\rm K}/\Q_3}(\alpha) \not\equiv 1({\rm mod} \ 3^6)$.
Now ${\rm N}_{{\rm K}/\Q_3}(\alpha) = \alpha^{[\rm{K}:\Q_3]} = \alpha^{54}$.
Since ${\rm U}_{1,\Q_3} \stackrel{(\ )^2}{\longrightarrow} {\rm U}_{1,\Q_3}$ is an isomorphism which preserves all filtration levels, $\alpha^2 \in {\rm U}_{2,\Q_3} \backslash {\rm U}_{3,\Q_3}$. By using Prop.\ref{p-power},
$${\rm U}_{2,\Q_3}/{\rm U}_{3,\Q_3}  \stackrel{(\ )^{27}}{\longrightarrow} {\rm U}_{5,\Q_3}/{\rm U}_{6,\Q_3}$$ is an isomorphism. This implies that ${\rm N}_{{\rm K}/\Q_3}(\alpha) \not\equiv 1({\rm mod} \ 3^6)$. \\

\noindent We have seen that, if ${\rm K} = \Q_p(\zeta_{p^m})$, ${\rm N}_{{\rm K}/\Q_p}({\rm U}_{h(m)+1, \rm K}) = {\rm N}_{{\rm K}/\Q_p}({\rm U}_{p^m,\rm K})  \subset {\rm U}_{m+1, \Q_p}$. By the property of Hasse-Herbrand function,(see sec.~\ref{results}) ${\rm NU}_{h(i),\rm L} \not\subset {\rm U}_{i+1,\rm F}$. We now deduce that ${\rm N}_{{\rm K}/\Q_p}({\rm U}_{h(m),{\rm K}}) = {\rm N}_{{\rm K}/\Q_p}({\rm U}_{p^m-1,{\rm K}}) \not\subset {\rm U}_{m+1, \Q_p}$. Thus proving that the congruence of the Theorem \ref{main} may not hold for the units which are not $p$-primary. 
\end{ex}


\noindent Supriya Pisolkar\\
\noindent Harish-Chandra Research Institute\\
\noindent Chhatnag Road, Jhunsi\\
\noindent Allahabad -211019, India.\\
\noindent {\it Email - supriya@mri.ernet.in}

\end{document}